\documentclass[11pt,leqno]{article}
 \title{On Furtw\"angler's theorems and second case of Fermat's Last Theorem} 
\author{Roland Qu\^eme}
\usepackage{amsmath,amsthm,amstext,amsbsy,eufrak}
\newtheorem{thm}{Theorem}[section]
\newtheorem{cor}[thm]{Corollary}

\newtheorem{lem}[thm]{Lemma}
\newtheorem{conj}{Conjecture}
\newtheorem{defi}{Definition}
\newtheorem{notas}{Notations}
\newtheorem{rema}{Remark}
\font\mathbb=msbm10
\newcommand{\N}{\mbox{\mathbb N}}

\newcommand{\Q}{\mbox{\mathbb Q}}
\newcommand{\Z}{\mbox{\mathbb Z}}
\newcommand{\be}{\begin{equation}}
\newcommand{\ee}{\end{equation}}
\newcommand{\bd}{\begin{displaymath}}
\newcommand{\ed}{\end{displaymath}}
\newcommand{\bn}{\begin{enumerate}}
\newcommand{\en}{\end{enumerate}}

\newcommand{\mk}{\mathfrak}
\newcommand{\ml}{\mathcal}
\newcommand{\mf}{\mathbf}

\newcommand{\ov}{\bar}
\newcommand{\bi }{\begin{itemize}}
\newcommand{\ei}{\end{itemize}}
\date{2013 Apr 23}
\begin{document}

\maketitle
%\tableofcontents
\maketitle
\begin{abstract}
This article, complement to the article [Que], deals with some  generalizations of  Futw\"angler's theorems for the second case of Fermat's Last Theorem (FLT2).
Let $p$ be an odd prime, $\zeta$ a $p$th primitive root of unity,  $K:=\Q(\zeta)$ and  $C\ell_K$   the class group of $K$.
A prime $q$ is said {\it  $p$-principal}
if the class $c\ell_K (\mk q_K)\in C\ell_K$  of any prime ideal 
$\mk q_K$ of $\Z_K$ over $q$   is the $p$th power of a class.
Assume that FLT2 fails for $(p,x,y,z)$ where $x, y, z$ are mutually coprime integers, $p$ divides $y$ and 
$x^p+y^p+z^p=0$.

Let $q$ be a prime dividing $\frac{(x^p+y^p)(y^p+z^p)(z^p+x^p)}{(x+y)(y+z)(z+x)}$ and  
 $\mk q_K$ be  any  prime ideal of $K$ over $q$.
We obtain  the $p$-power residue symbols  relations
$$\Big(\frac{p}{\mk q_K}\Big)_{\!\! K}=\Big(\frac{1-\zeta^j}{\mk q_K}\Big)_{\!\! K} \mbox{\ for\ } j=1,\dots,p-1.$$
As an application, 
we prove that: if Vandiver's conjecture holds  for $p$ then  
$q$ is a $p$-principal prime. 

Similarly, let $q$ be a prime dividing dividing $\frac{(x^p-y^p)(y^p-z^p)(z^p-x^p)}{(x-y)(y-z)(z-x)}$ and  
$\mk q_K$ be the  prime ideal of $K$ over $q$ dividing $(x\zeta-y)(z\zeta-y)(x\zeta -z)$.
We give an explicit formula for  the $p$-power residue symbols 
$\big(\frac{\epsilon_{k}}{\mk q_K}\big)_{\!\! K}$ for  all $k$ with $1<k\leq\frac{p-1}{2},$
where 
$\epsilon_k$ is the cyclotomic unit given by $\epsilon_k=:\zeta^{(1-k)/2}\cdot\frac{1+\zeta^k}{1+\zeta}.$

The principle of proofs rely on the $p$-Hilbert class field theory.
\end{abstract}

\footnote{

- {\it Keywords:} Fermat's Last Theorem, cyclotomic fields, cyclotomic
units, 
class field theory, Vandiver's and Furtw\"angler's theorems

- {\it AMS subject classes :} 11D41, 11R18, 11R37}

\maketitle
\section{  Introduction } \label{s1012111}

\subsection{General  notations and definitions } \label{definitions}

\bi
\item
Let  $p>3$ be a prime, $\zeta:=e^{\frac{2\pi i}{p}}$, $K:=\Q(\zeta)$
the  $p$th cyclotomic number field, $\Z_K$ the ring of integers  of $K$, 
 and  $\mk p=(1-\zeta)\Z_K$  the prime ideal of $\Z_K$ over $p$.
Let $g:={\rm Gal}(K/\Q)$, for $k \not\equiv 0 \mod p$ and  $s_k :
\zeta \rightarrow\zeta^k$   
 the $p-1$ distinct elements of $g$. 
\item
Let $C\ell_K$, $C\ell$ and $C\ell^-$   be respectively the class group of $K$, 
the $p$-class group of $K$ and  the negative part of the $p$-class group of $K$.
For any ideal  $\mk a$  of $K$, 
let us note $c\ell_K(\mk a)$, $c\ell(\mk a), c\ell^-(\mk a)$ be respectively the class of $\mk a$ in $C\ell_K$,  $C\ell$ and $C\ell^-$.
%\item
%A number $\alpha\in K^\times$ prime to $p$, such that $\alpha\Z_K$ is
%the $p$th power of an ideal,  is called a pseudo-unit.
%The pseudo-unit $\alpha$ is   $p$-primary  (i.e. 
%the extension $K(\sqrt[p]{\alpha})/K$ is unramified at $\mk p$) 
%if and only if $\alpha$ is congruent 
%to a $p$-power $\mod \mk p^p$, see [Gr2] lem 2.1.
\item
A prime $q$ is said {\it  $p$-principal}
if the class $c\ell_K (\mk q_K)\in C\ell_K$  of any prime ideal 
$\mk q_K$ of $\Z_K$ above  $q$ is the $p$th power of a class, which is 
equivalent to $\mk q_K = {\mk a}^p (\alpha)$, for an ideal ${\mk a}$
of $K$ and an $\alpha \in K^\times$.
This contains the case where  the class $c\ell_K (\mk q_K)$ is of order
coprime with  $p$.
\item For any $\alpha\in K$ and prime ideal $\mk q_K$ of $K$, we use the $p$th power residue symbol notation 
$\big(\frac{\alpha}{{\mathfrak q}_K}\big)_{\!\!K}$. 
\item 
We will adopt  in the sequel the following notations for  an hypothetic
counterexample to $FLT2$. We say that $FLT2$  would fail for $(p,x,y,z)$ if we had 
$$x^p+y^p+z^p=0,$$
 with $x,y,z\in\Z\backslash\{0\}$ pairwise coprime and $p$ dividing  $y$.
\ei
\subsection{Main results}
Let $q$ be a prime dividing $\frac{(x^p+y^p)(y^p+z^p)(z^p+x^p)}{(x+y)(y+z)(z+x)}$ and  
 $\mk q_K$ be  any  prime ideal of $K$ over $q$.
We obtain  the $p$-power residue symbols  relations (see theorem \ref{t1107231})
$$\Big(\frac{p}{\mk q_K}\Big)_{\!\! K}=\Big(\frac{1-\zeta^j}{\mk q_K}\Big)_{\!\! K} \mbox{\ for\ } j=1,\dots,p-1.$$
%$\big(\frac{x\zeta+y}{\mk q_K}\big)_{\!\!K}$,
%$\big(\frac{z\zeta +y}{{\mathfrak q}_K}\big)_{\!\!K}$
%and  $\big(\frac{x\zeta+z}{{\mathfrak q}_K}\big)_{\!\!K}$ 
%depending only on the cyclotomic units 
%$\varpi_k:=\zeta^{(1-k)/2}\cdot\frac{1-\zeta^k}{1-\zeta},\ 1 < k<\frac{p-1}{2}.$ 
As an application, 
we prove that: if Vandiver's conjecture fails for $p$ then  
$q$ is a $p$-principal prime (see theorem \ref{t1111221}).
% and $\mk q_K$  splits  totally in the $p$-Kummer extension
%$K\big (\sqrt[p] {<\varpi_k>_{k=2,\dots,(p-3)/2}}\ \big).$

Similarly, let $q$ be a prime dividing dividing $\frac{(x^p-y^p)(y^p-z^p)(z^p-x^p)}{(x-y)(y-z)(z-x)}$ and  
$\mk q_K$ be the  prime ideal of $K$ over $q$ dividing $(x\zeta-y)(z\zeta-y)(x\zeta -z)$.
We give an explicit formula for the $p$-power residue symbols 
$\big(\frac{\epsilon_{k}}{\mk q_K}\big)_{\!\! K}$ for  all $k$ with $1<k\leq\frac{p-1}{2},$
where 
$\epsilon_k$ is the cyclotomic unit given by $\epsilon_k=:\zeta^{(1-k)/2}\cdot\frac{1+\zeta^k}{1+\zeta}$ (see theorem \ref{t1107232}).

%***
%Let $q$ be a prime dividing dividing $\frac{(x^p+y^p)(y^p+z^p)(z^p+x^p)}{(x+y)(y+z)(z+x)}$ and  
% $\mk q_K$ a prime ideal of $K$ over $q$.
%We give an explicit computation of the $p$-power residue symbols $\big(\frac{x\zeta+y}{\mk q_K}\big)_{\!\!K}$,
%$\big(\frac{z\zeta +y}{{\mathfrak q}_K}\big)_{\!\!K}$
%and  $\big(\frac{x\zeta+z}{{\mathfrak q}_K}\big)_{\!\!K}$ 
%depending only on the cyclotomic units 
%$\varpi_k:=\zeta^{(1-k)/2}\cdot\frac{1-\zeta^k}{1-\zeta},\ 1 < k<\frac{p-1}{2}.$ As an application we
%prove that: {\it  if Vandiver's conjecture fails for $p$ then 
%$q$ is a $p$-principal  prime and $\mk q_K$  splits  totally in the $p$-Kummer extension
%$K\big (\sqrt[p] {<\varpi_k>_{k=2,\dots,(p-3)/2}}\ \big).$ }
%Similarly, let $q$ be a prime dividing  $\frac{(x^p-y^p)(y^p-z^p)(z^p-x^p)}{(x-y)(y-z)(z-x)}$ and  
%$\mk q_K$ a prime ideal of $K$ over $q$.
%We give an explicit computation of the $p$-power residue symbols $\big(\frac{x\zeta-y}{\mk q_K}\big)_{\!\!K}$,
%$\big(\frac{z\zeta -y}{{\mathfrak q}_K}\big)_{\!\!K}$
%and  $\big(\frac{x\zeta-z}{{\mathfrak q}_K}\big)_{\!\!K}$
%depending only on the cyclotomic units 
%$\zeta^{(1-k)/2}\cdot\frac{1+\zeta^k}{1+\zeta},\ 1 <k < \frac{p-1}{2}.$
%***
This article is a complement to the article [GQ] dealing
with {\it Strong Fermat's Last Theorem conjecture (SFLT)} and article  [Que] dealing 
with {\it second case of Strong Fermat's Last Theorem conjecture (SFLT2)}.

\section {Detailed results and proofs}\label{s1011101}
\label{FLT2}
%This section details    the particular case $FLT2$ of $SFLT2$.  
%The last subsection  focus more particularly on some strong properties of
%the primes
%dividing $\frac{(x^p+y^p)(y^p+z^p)(z^p+x^p)}{(x+y)(y+z)(z+x)}$ 
%or  $\frac{(x^p-y^p)(y^p-z^p)(z^p-x^p)}{(x-y)(y-z)(z-x)}$
%if FLT2 failed for $(p,x,y,z)$ with $p|y$.
We give at first a general lemma.
\begin{lem}\label{l3}
Suppose that $FLT2$ fails for $(p,x,y,z)$ with $p|y$.  If  $q\not=p$ satisfies  
$$y\equiv 0\mod q \mbox{\ and\ }
 x+z\not\equiv 0\mod q,$$ then $q-1\equiv 0\mod  p^2$.
\begin{proof}$ $
\bi
\item
  From Barlow-Abel relations 
$$x+z =p^{\nu p-1} y_0^p,\ \frac{x^p+z^p}{x+z}= py_1^p,\ y=- p^\nu 
y_0y_1, \ \nu\geq 1,$$
\item
Suppose that  $q|\frac{x^p+z^p}{x+z}$ with $p$ prime to $\kappa$  and
search for a contradiction:
let $\mk q_K$  be a prime ideal of $\Z_K$ lying over $q$. From $q|y$ 
and the Barlow-Abel relation  $x+y=z_0^p$, we have
$$\Big(\frac{x}{\mk  q_K}\Big)_{\!\! K}=\Big(\frac{x+y}{\mk
q_K}\Big)_{\!\! K}=
\Big(\frac{z_0^p}{\mk  q_K}\Big)_{\!\! K}=1.$$ 
Similarly $\Big(\frac{z}{\mk q_K}\Big)_{\!\!K}=1$, 
so $x^{(q-1)/p} - z^{(q-1)/p}\equiv 0\mod  \mk q_K$. 
We get $$q\ |\  x^{(q-1)/p}-z^{(q-1)/p} \mbox{\  and \ } 
q\ |\ x^p+z^p.$$
\item  
If we suppose $\kappa=\frac{q-1}{p}$ prime 
to $ p$, we have $\kappa=\frac{q-1}{p}$ even  
and $ x^\kappa \equiv (-z)^\kappa \mod  q$  and $ x^p \equiv  (-z)^p
\mod   q$,  
thus  $q \ |\  x+z$  by a  B\'ezout relation between  $p$ and $n$
(absurd). 
\ei
\end{proof}
\end{lem}

\subsection{On the   primes $q$ dividing 
$\frac{(x^p+y^p)(y^p+z^p)(z^p+x^p)}{(x+y)(y+z)(z+x)}$} $ $
\bn
\item
We assume that $FLT2$ fails for $(p,x,y,z)$.
This section contains
some general  strong properties  of the primes $q$ dividing  
$\frac{(x^p+y^p)(y^p+z^p)(z^p+x^p)}{(x+y)(y+z)(z+x)}$ complementary to Furtw\"angler's theorems.
Here, we don't assume  that $q$ is $p$-principal or not, 
thus this subsection brings complementary informations to corollary
2.7 of [Que]. 
\item
Let us define the  totally real cyclotomic units   
$$\varpi_{a}=:\zeta^{(1-a)/2}\cdot\frac{1-\zeta^{a}}{1-\zeta},\ 1 \leq
a\leq p-1,$$
where this definition implies  $\varpi_1=1$.
Recall that the cyclotomic units of $K$ are generated by the $\varpi_a$
for $1<a<\frac{p}{2}$.
We have $\varpi_a=-\varpi_{p-a}$: indeed we have
$\varpi_a=\zeta^{(1-a)/2}\cdot\frac{1-\zeta^a}{1-\zeta}$ 
and $\varpi_{p-a}=
\zeta^{(1-(p-a))/2}\cdot\frac{1-\zeta^{p-a}}{1-\zeta}=\zeta^{(1+a)/2}\cdot\frac{1-\zeta^{-a}}{1-\zeta}
=\zeta^{1-a)/2}\cdot\frac{\zeta^a-1}{1-\zeta}=-\varpi_a$.
\en

\begin{lem}\label{l1108221}
Assume that FLT2 fails for $(p,x,y,z)$ with $p|y$ . 
Let $\mk q_K$ be a prime  ideal of $\Z_K$ such that $x\zeta+y\equiv 0 \mod
\mk q_K$ (or $z\zeta+  y\equiv 0\mod \mk q_K$).
Then
$$q\equiv 1\mod p^2 \mbox{\ and\ } \Big(\frac{\zeta}{\mk q_K}\Big)_{\!\!
K}=\Big(\frac{p}{\mk q_K}\Big)_{\!\! K}
=\Big(\frac{1-\zeta}{\mk q_K}\Big)_{\!\! K}=1.$$

\begin{proof}$ $
\bi
\item
 Suppose that $x\zeta +y\equiv 0\mod \mk q_K$.
We have $q|z$, so $q\equiv 1\mod p^2$ from First Furtw\"angler's theorem, so
$\Big(\frac{\zeta}{\mk q_K}\Big)_{\!\! K}=1$ and  
$\Big(\frac{x}{\mk q_K}\Big)_{\!\! K}=\Big(\frac{y}{\mk q_K}\Big)_{\!\!
K},$ so
$\Big(\frac{x+z}{\mk q_K}\Big)_{\!\! K}=\Big(\frac{y+z}{\mk
q_K}\Big)_{\!\! K},$   
so 
$$\Big(\frac{p^{\nu p-1} y_0^p}{\mk q_K}\Big)_{\!\!
K}=\Big(\frac{x_0^p}{\mk q_K}\Big)_{\!\! K}\mbox{\ with \ }\nu\in \N_{\geq 1},$$ 
from Barlow-Abel relations, 
and finally 
$\Big(\frac{p}{\mk q_K}\Big)_{\!\! K}=1.$
In the other hand, we have 
$$x+y=z_0^p\equiv x(1-\zeta)\equiv (x+z)(1-\zeta)\equiv p^{\nu
p-1}y_0^p(1-\zeta)\mod \mk q_K,$$
so
$$\Big(\frac{1-\zeta}{\mk q_K}\Big)_{\!\! K}=1.$$
\item
 Suppose  that $z\zeta+y\equiv 0\mod \mk q_K$.  The proof  is similar
with $z$ in place of $x$.
\ei
\end{proof}
\end{lem}

\begin{lem}\label{l1304113}
Suppose that $FLT2$ fails for $(p,x,y,z)$ with $p|y$ .
Let $q\not= p$ be a prime and  $\mk q_K$ be a  prime ideal of $\Z_K$ over
$q$. Then we have 
for $k=1,\dots,p-2$: 
\bn
\item
 If $\mk q_K$ divides $x\zeta+y$  then $\Big(\frac{x+\zeta^k y}{\mk
q_K}\Big)_{\!\! K}
=\Big(\frac{\varpi_{k+1}}{\mk q_K}\Big)_{\!\! K}$. 
\item
If $ \mk q_K$ divides $z\zeta+y$  then $\Big(\frac{z+\zeta^k y}{\mk
q_K}\Big)_{\!\! K}
=\Big(\frac{\varpi_{k+1}}{\mk q_K}\Big)_{\!\! K}$.
\item
 If $\mk q_K$ divides $x\zeta+z$ and $p\ |\ y$   
then $\Big(\frac{x+\zeta^k z}{\mk q_K}\Big)_{\!\! K}\Big(\frac{p }{\mk
q_K}\Big)_{\!\! K}
=\Big(\frac{\varpi_{k+1}}{\mk q_K}\Big)_{\!\! K}$.
\en
\begin{proof}$ $
\bn
\item
 From  $x\zeta + y\equiv 0\mod \mk q_K$ we get 
 
$$x+\zeta^k y\equiv x(1-\zeta^{k+1})\mod \mk q_K,\ k=1,\dots,p-2.$$
thus 
$$\frac{x+\zeta^k y}{x+y}\equiv \frac{1-\zeta^{k+1}}{1-\zeta}\mod \mk q_K,
\mbox{\ for\ } k=1,\dots,p-2.$$
In the other hand,
$\varpi_{k+1}=\zeta^{(1-(k+1))/2}\cdot\frac{1-\zeta^{k+1}}{1-\zeta}$ 
is a totally real cyclotomic unit, 
so  
$$\frac{x+\zeta^k y}{x+y}\equiv \varpi_{k+1}\zeta^{k/2}\mod\mk q_K,\mbox{\
for\ } k=1,\dots p-2,$$ 
so 
$$\Big (\frac{x+\zeta^k y}{\mk q_K}\Big)_{\!\! K}
=\Big(\frac{\varpi_{k+1}}{\mk q_K}\Big)_{\!\!
K}\Big(\frac{\zeta^{k/2}}{\mk q_K}\Big)_{\!\! K} 
\mbox{\ for\ } k=1,\dots,p-2,$$
because $x+y\in K^{\times p}$ and finally
$$\Big (\frac{x+\zeta^k y}{\mk q_K}\Big)_{\!\! K}
=\Big(\frac{\varpi_{k+1}}{\mk q_K}\Big)_{\!\! K}\mbox{\ for\ }
k=1,\dots,p-2,$$
because $q\equiv 1\mod p^2$ obtained by the first Theorem of Furtw\"angler.
\item
 The proof is similar to  {\it item 1.} with $z$ in place of $x$. 
\item
 In that case we have $x+z=p^{\nu p-1}y_0^p$ with $\nu>0$ 
and so $x+z\in p^{-1}K^{\times p}$ and  $p^2|q-1$  as proved 
in   lemma \ref{l3}.
\en
\end{proof}
\end{lem}
%\begin{rema}
%{\rm This  property of the primes $q$ dividing $\frac{x^p+y^p}{x+y}$  
%(or  $\frac{z^p+y^p}{z+y}$, or  $\frac{x^p+z^p}{x+z})$   is strong
%because $x\zeta+ y$ (or $z\zeta+ y$,
% or $\frac{x+\zeta z}{1-\zeta})$),  
%and $\varpi_{k+1}$ of $\Z_K$  
%are  pseudo-units not linearly connected in the action of $\F_p[g]$.
%}
%\end{rema}

\begin{thm}\label{t1107231}
Assume   that   the second case of FLT fails for $(p,x,y,z)$ with $p|y$.
Let $q$ be a prime dividing $\frac{x^p+y^p}{x+y}$ (or
$\frac{z^p+y^p}{z+y}$ 
or $\frac{x^p+z^p}{x+z}$).
Let $\mk q_K$ be \underline{the} prime ideal of $\Z_K$ over $q$
dividing $x\zeta+ y$ (or $z\zeta+y$ or $x\zeta+z$).
 
If the $p$-class  $c\ell(\mk q_K)\in C\ell^-$  we have: 
\bn
\item \label{i1}
The prime $q$ satisfies the congruence $q\equiv 1\mod p^2$.
\item
$\mk q_K$ satisfies  the  following power residue  symbols  values:
\bn
\item \label{i2a}
If $\mk q_K|x\zeta+y$ (or $z\zeta+y$) then $$\Big(\frac{p}{\mk q_K}\Big)_{\!\! K}=\Big(\frac{1-\zeta^j}{\mk q_K}\Big)_{\!\! K}=1\mbox{\ for\ } j=1,\dots,p-1.$$
\item \label{i2b}
If $\mk q_K|x\zeta+z$  then $$\Big(\frac{p}{\mk q_K}\Big)_{\!\! K}=\Big(\frac{1-\zeta^j}{\mk q_K}\Big)_{\!\! K} \mbox{\ for\ } j=1,\dots,p-1.$$
\item \label{i2c}
If Vandiver's conjecture holds for $p$, the prime $q$ is $p$-principal.
\en 
\en
\begin{proof}$ $
\begin{itemize}
  \item 
If  $q|\frac{x^p+y^p}{x+y}\frac{x^p+y^p}{x+y}$, from Furtwangler's
First theorem, we get  $q\equiv 1\mod p^2$.
We derive that $\Big(\frac{\zeta}{\mk q_K}\Big)_{\!\! K}=1$ and from lemma
\ref{l1108221}  that  $\Big(\frac{p}{\mk q_K}\Big)_{\!\! K}=1$.
If $q|\frac{x^p+z^p}{x+z}$ then, $q\equiv 1\mod p^2$  from  lemma \ref{l3}, which proves {\it item  \ref{i1}} of the statement. 
\item  Suppose $q|\frac{x^p+y^p}{x+y}$.
\begin{itemize}
\item 
From   previous lemma \ref{l1304113},  we have 
$$\Big (\frac{x+\zeta^k y}{\mk q_K}\Big)_{\!\! K}
=\Big(\frac{\varpi_{k+1}}{\mk q_K}\Big)_{\!\! K} 
\mbox{\ for\ } k=1,\dots,p-2,$$
and also, with $p-k$ in place of $k$, 
$$\Big (\frac{x+\zeta^{p-k} y}{\mk q_K}\Big)_{\!\! K}
=\Big(\frac{\varpi_{p-k+1}}{\mk q_K}\Big)_{\!\! K} 
\mbox{\ for\ } p-k=1,\dots,p-2,$$
so
\be\label{e1304111}
\Big (\frac{\frac{x+\zeta^ ky}{x+\zeta^{p-k}y} }{\mk q_K}\Big)_{\!\! K}
=\Big(\frac{\varpi_{k+1}\varpi_{p-k+1}^{-1}}{\mk q_K}\Big)_{\!\! K} 
\mbox{\ for\ } p-k=1,\dots,p-2,\ee
\item
For $2\leq k\leq p-2$, we can write 
$$x+\zeta^k y= A_k B_k\alpha^p,$$ with $\alpha\in K^{\times p}$,  
pseudo-units $A_k, B_k$ verifying $A_k^{s_{-1}+1}\in K^{\times p}$ 
and $B_k^{s_{-1}-1}\in K^{\times p}$ where we recall that $s_k$ is the
$\Q$-isomorphism $s_k: \zeta\rightarrow \zeta^k$ of $K$. 
Let $\Big(\frac{A_k}{\mk q_K}\Big)_{\!\! K}=\zeta^w$, we get 
$$\Big(\frac{A_k^{s_{-1}}}{s_{-1}(\mk q_K)}\Big)_{\!\! K}
=\Big(\frac{A_k^{-1}}{s_{-1}(\mk q_K)}\Big)_{\!\! K}=\zeta^{-w},$$
so $$\Big(\frac{A_k}{s_{-1}(\mk q_K)}\Big)_{\!\! K}=\zeta^{w},$$ 
and so $\Big(\frac{A_k}{\mk q_K s_{-1}(\mk q_K)}\Big)_{\!\! K}=\zeta^{2w}$.
But $c\ell(\mk q_K)\in C\ell^-$, so $(\mk q_Ks_{-1}(\mk
q_K))^n\Z_K=\beta\Z_K$ 
with $\beta\in \Z_K$ and a certain integer $n$  coprime with $p$. Then 
$$\Big(\frac{A_k}{\mk q_K^n s_{-1}(\mk q_K)^n}\Big)_{\!\! K}=
\Big(\frac{A_k}{\beta}\Big)_{\!\! K} =1,$$ 
because $A_k$ is a $p$-primary pseudo-unit 
(for instance by application of Artin-Hasse reciprocity law), so $w=0$ and 
$\Big(\frac{A_k}{\mk q_K}\Big)_{\!\! K}=1$.
\item  We get $\frac{x+\zeta^k y}{x+\zeta^{p-k} y}\in A_k^2\times K^{\times
p}$, so
\be\label{e1111105} 
\Big(\frac{x+\zeta^k y}{\mk q_K}\Big)_{\!\!}
=\Big(\frac{x+\zeta^{p-k}y}{\mk q_K }\Big)_{\!\! K}\mbox{\ for \ }k=2,\dots,p-2.
\ee
%which leads  from (\ref{e1111101}) and (\ref{e1111102}) to
%\be\label{e1107081}
%\Big(\frac{\epsilon_{p-k+1}}{\mk q_K}\Big)_{\!\! K}
%=\Big(\frac{\zeta^{k}}{\mk q_K}\Big)_{\!\!
%K}\Big(\frac{\epsilon_{k+1}}{\mk q_K}\Big)_{\!\! K} 
%\mbox{\ for\ } k=2,\dots,p-2.
%\ee
%The relation (\ref{e1111105}) of theorem \ref{t1107232} can be obtained with exactly the same proof
%\be\label{e1111107} 
%\Big(\frac{x+\zeta^k y}{\mk q_K}\Big)_{\!\!}
%=\Big(\frac{x+\zeta^{p-k}y}{\mk q_K }\Big)_{\!\! K}\mbox{\ for \ }k=2,\dots,p-2,
%\ee
which leads to
$$\Big(\frac{\varpi_{k+1}}{\mk q_K}\Big)_{\!\! K}
=\Big(\frac{\varpi_{p-k+1}}{\mk q_K}\Big)_{\!\! K}  \mbox{\ for\ }
k=2,\dots,p-2.$$
\item
We have  seen above that  $\varpi_{k+1}= -\varpi_{p-k-1}$
so 
$$\Big(\frac{\varpi_{k+1}}{\mk q_K}\Big)_{\!\! K}=
\Big(\frac{\varpi_{p-k-1}}{\mk q_K}
\Big)_{\!\! K}\mbox{\ for\ }k=2,\dots,p-2.$$
Then, gathering these relations involving the units $\varpi_{k+1}, \varpi_{p-k-1},\varpi_{p-k+1}$,  we get 
$$\Big(\frac{\varpi_{p-k+1}}{\mk q_K}\Big)_{\!\! K}
=\Big(\frac{\varpi_{p-k-1}}{\mk q_K}\Big)_{\!\! K}  \mbox{\ for\ }
k=2,\dots,p-2.$$

\item Starting from $k=2$ we get  for $k=2,4,\dots,p-3,$
$$\Big(\frac{\varpi_{p-1}}{\mk q_K}\Big)_{\!\! K}
=\Big(\frac{\varpi_{p-3}}{\mk q_K}\Big)_{\!\! K} =\dots =
\Big(\frac{\varpi_{2}}{\mk q_K}\Big)_{\!\! K} =1,$$
because we get directly $\Big(\frac{\varpi_{p-1}}{\mk q_K}\Big)_{\!\!
K}=1$ from its definition.
Starting from $k=3$ we get for $k=3,5,\dots, p-2,$
$$\Big(\frac{\varpi_{p-2}}{\mk q_K}\Big)_{\!\! K}
=\Big(\frac{\varpi_{p-4}}{\mk q_K}\Big)_{\!\! K} =\dots =
\Big(\frac{\varpi_{1}}{\mk q_K}\Big)_{\!\! K} =1,$$
because we get directly $\Big(\frac{\varpi_1}{\mk q_K}\Big)_{\!\! K}=1$
from its definition. 
Therefore we get $$\Big(\frac{\varpi_i}{\mk q_K}\Big)_{\!\! K}=1\mbox{\
for\ }i=1,\dots, p-1.$$ 
So,  we get $$\Big(\frac{1-\zeta^i}{\mk q_K}\Big)_{\!\! K}=\Big(\frac{1-\zeta}{\mk q_K}\Big)_{\!\! K}\mbox{\
for\ }i=1,\dots, p-1.$$ 
and finally
we find again $\Big(\frac{p}{\mk q_K}\Big)_{\!\!
K}=\Big(\frac{1-\zeta}{\mk q_K}\Big)_{\!\! K},$ seen in lemma 
\ref{l1108221}. 

From  lemma \ref{l1108221} we have also  $\Big(\frac{1-\zeta}{\mk
q_K}\Big)_{\!\! K}=1$ if $\mk q_K|x\zeta+y$ (or $\mk q_K|z\zeta+y$),
which proves {\it item 2.a}  for $q|\frac{(x^p+y^p)(z^p+y^p}{(x+y)(z+y)}$.

\item If Vandiver's conjecture holds for $p$ the $p$-primary units corresponding to $C\ell^-$  are all generated by the 
$\varpi_i,\ \ i=1,\dots,\frac{p-1}{2}$.
Therefore,  the result 
$\Big(\frac{\varpi_{i}}{\mk q_K}\Big)_{\!\! K}=1\mbox{\ for\
}i=1,\dots,p-1$ obtained 
and the assumption that $c\ell(\mk q_K)\in C\ell^-$ 
imply  that $\mk q_K$ is $p$-principal  
(application of  the  decomposition  and reflection theorems in the $p$-Hilbert class
field of $K$),
if not it should be possible to find  some integers $n_1,\dots,n_{(p-3)/2}\not\equiv 0\mod p,$
such that the   $p$-primary unit $\varpi=\prod_{i=1}^{(p-3)/2}\varpi_i^{n_i}$ verifies 
$\Big(\frac{\varpi}{\mk q_K}\Big)_{\!\! K}\not=1$, contradiction
which proves {\it item 2.c}  for $q|\frac{(x^p+y^p)(z^p+y^p}{(x+y)(z+y)}$.
\end{itemize}

\item Suppose at last  that $q|\frac{x^p+z^p}{x+z}$:
 If $\mk q_K|x\zeta+z$ and $p\ |\ y$   then 
$$\Big(\frac{x+\zeta^k z}{\mk q_K}\Big)_{\!\! K}\Big(\frac{p }{\mk
q_K}\Big)_{\!\! K}
=\Big(\frac{\varpi_{k+1}}{\mk q_K}\Big)_{\!\! K},$$ 
(seen in lemma \ref{l1304113} {\it item 3.}) and similarly 
$$\Big(\frac{x+\zeta^{p-k} z}{\mk q_K}\Big)_{\!\! K}\Big(\frac{p }{\mk
q_K}\Big)_{\!\! K}
=\Big(\frac{\varpi_{p-k+1}}{\mk q_K}\Big)_{\!\! K},$$ 
so  we get again the relation (\ref{e1304111})
$$\Big(\frac{\frac{x+\zeta^k z}{x+\zeta^{p-k} z}}{\mk q_K}\Big)_{\!\! K}
=\Big(\frac{\varpi_{k+1}\varpi_{p-k+1}^{-1}}{\mk q_K}\Big)_{\!\! K}.$$
In the other hand $\frac{x+\zeta^k z}{x+\zeta^{p-k} z}= \zeta^k A$  
where $A$ is also a  $p$-primary pseudo unit with $A^{s_{-1}+1}\in K^{\times p}$.
Then the end of the proof is similar to the previous cases 
 $q|\frac{(x^p+y^p)(z^p+y^p)}{(x+y)(z+y)}$
taking into
account that we know that $p^2|q-1$, so $\Big(\frac{\zeta^k}{\mk q_K}\Big)_{\!\! K}=1$,
which proves {\it  items   2b. and 2c.} of the statement if $q|\frac{x^p+z^p}{x+z}$.
\end{itemize}
\end{proof}
\end{thm}

\begin{rema}
{\rm 
In the case of an hypothetic solution $(x,y,z),\ p|y$ of the FLT2
equation, 
for the primes $q$ with  $c\ell(\mk q_K)\in C\ell^-$ and $\mk
q_K|x\zeta+y$ (or $z\zeta+y)$,  the theorem \ref{t1107231} can be
considered as  a   reciprocal  statement to corollary 2.7 of [Que]
in which $(u,v)=(x,y)$ or $(z,y)$ for $x,y,z, \ p|y$  hypothetic  solution
of the Fermat's equation.
In particular,  we have proved:}
\end{rema}
\begin{thm}\label{t1111221}
Assume that  Vandiver's conjecture holds for $p$ and that  the second case of FLT fails for $(p,x,y,z)$.
Then all the primes $q\not=p$ dividing $\frac{(x^p+y^p)(y^p+z^p)(z^p+x^p)}{(x+y)(y+z)(z+x)}$
are $p$-principal.
\end{thm}

\subsection{Some properties of  the  primes $q$ 
dividing $\frac{(x^p-y^p)(y^p-z^p)(z^p-x^p)}{(x-y)(y-z)(z-x)}$ }
\bn
\item
We assume that the second case  $FLT2$ fails for $(p,x,y,z)$ with $p|y$.
This subsection contains
some general    properties  of decomposition of  the primes $q$
dividing  
$\frac{(x^p-y^p)(y^p-z^p)(z^p-x^p)}{(x-y)(y-z)(z-x)}$ in certain
$p$-Kummer extensions.
Here, we don't assume  that $q$ is $p$-principal or not, 
thus this subsection brings complementary informations to $SFLT2$  corollary
2.5 in [Que]. 
Note that, here,  Furtw\"angler's theorems cannot be applied  to these primes $q$, so  we
cannot assume that $p^2$ divides $q-1$.
\item
Let us define   the  totally real cyclotomic units   
$$\epsilon_{a}=:\zeta^{(1-a)/2}\cdot\frac{1+\zeta^{a}}{1+\zeta},\ 1 \leq
a\leq p-1,$$
where we note that $\epsilon_1=1$ 
and that
\be\label{e1111103} 
\varepsilon_{p-a}=\zeta^{(1-(p-a))/2}\cdot\frac{1+\zeta^{p-a}}{1+\zeta}= \zeta^{(1+a)/2}\cdot\frac{1+\zeta^{-a}}{1+\zeta}= \zeta^{(1-a)/2}\frac{1+\zeta^a}{1+\zeta}=\varepsilon_a.
\ee
\en
\begin{lem}\label{l1304111}
Suppose that $FLT2$ fails for $(p,x,y,z)$ with $p|y$ .
Let $q\not= p$ be a prime and  $\mk q_K$ be a  prime ideal of $\Z_K$ over
$q$. Then we have 
for $k=1,\dots,p-1$: 
\bn 
\item
If $\mk q_K|x\zeta-y$  then $\Big(\frac{x+\zeta^k y}{\mk
q_K}\Big)_{\!\! K}
=\Big(\frac{\zeta^{k/2}}{\mk q_K}\Big)_{\!\!
K}\Big(\frac{\epsilon_{k+1}}{\mk q_K}\Big)_{\!\! K}$. 
\item
If $ \mk q_K|z\zeta-y$  then $\Big(\frac{z+\zeta^k y}{\mk
q_K}\Big)_{\!\! K}
=\Big(\frac{\zeta^{k/2}}{\mk q_K}\Big)_{\!\!
K}\Big(\frac{\epsilon_{k+1}}{\mk q_K}\Big)_{\!\! K}$.
\item
If $\mk q_K|x\zeta-z$    
then $\Big(\frac{x+\zeta^k z}{\mk q_K}\Big)_{\!\! K}\Big(\frac{p }{\mk
q_K}\Big)_{\!\! K}
=\Big(\frac{\zeta^{k/2}}{\mk q_K}\Big)_{\!\!
K}\Big(\frac{\epsilon_{k+1}}{\mk q_K}\Big)_{\!\! K}$.
\en
\begin{proof}$ $
\bn
\item
From  $x\zeta - y\equiv 0\mod\mk  q_K$ we get   
$$x+\zeta^k y\equiv x(1+\zeta^{k+1})\mod \mk q_K,\ k=1,\dots,p-1.$$
thus 
$$\frac{x+\zeta^k y}{x+y}\equiv 
\frac{1+\zeta^{k+1}}{1+\zeta}\mod \mk q_K,\mbox{\ for\ } k=1,\dots,p-1.$$
In the other hand, for $1\leq k\leq p-2$ then  
$\epsilon_{k+1}=\zeta^{(1-(k+1))/2}\cdot\frac{1+\zeta^{k+1}}{1+\zeta}$ 
is a totally real cyclotomic unit, so  
$\frac{x+\zeta^k y}{x+y}\equiv \epsilon_{k+1}\zeta^{k/2}\mod\mk q_K,\
k=1,\dots p-1$, 
and finally 
$$\Big (\frac{x+\zeta^k y}{\mk q_K}\Big)_{\!\! K}
=\Big(\frac{\zeta^{k/2}}{\mk q_K}\Big)_{\!\!
K}\Big(\frac{\epsilon_{k+1}}{\mk q_K}\Big)_{\!\! K} 
\mbox{\ for\ } k=1,\dots,p-2,$$
because  $x+y\in K^{\times p}$.
\footnote{ We don't know here if $p^2|q-1$.}
\item
The proof is similar with $z$ in place of $x$. 
\item
In that case we have $x+z=p^{\nu p-1}y_0^p$ with $\nu>0$ and so
$x+z\in p^{-1}K^{\times p}$.
\en
\end{proof}
\end{lem}

\begin{thm}\label{t1107232}
Suppose  that   the second case of FLT fails for $(p,x,y,z)$ with $p|y$.
Let $q$ be a prime dividing $\frac{x^p-y^p}{x-y}$ (or
$\frac{y^p-z^p}{y-z}$).
Let $\mk q_K$ be \underline{the} prime ideal of $\Z_K$ over $q$
dividing $x\zeta- y$ (or $z\zeta-y$). 
Assume that  the $p$-class   $c\ell(\mk q_K)\in C\ell^-$.
\footnote{As soon as Vandiver's conjecture is true for $p$, this assumption is verified.}
\bn
\item
If $p^2\not|\ q-1$   then  $q$ is non $p$-principal and satisfies 
$$\Big(\frac{\epsilon_{p-2k'-1}}{\mk q_K}\Big)_{\!\! K}
=\Big(\frac{\zeta^{-k'(k'+1)}}{\mk q_K}\Big)_{\!\! K}
\mbox{\ for\ }1\leq k'\leq \frac{p-3}{2},$$
and 
$$\Big(\frac{\epsilon_{p-2k'}}{\mk q_K}\Big)_{\!\! K}
=\Big(\frac{\zeta^{\frac{1}{4}-k'^2}}{\mk q_K}\Big)_{\!\! K}
\mbox {\ for\ }1\leq k'\leq \frac{p-3}{2}.$$
\item
If $p^2|\ q-1$ then $q$ satisfies
$$\Big(\frac{1+\zeta^j}{\mk q_K}\Big)_{\!\! K}=1\mbox {\ for\ } j=1,\dots
p-1.$$
\en

\begin{proof}$ $
\bn
\item
Let us  suppose at first that $p^2\not|\ q-1$: we know that $q$ is
non $p$-principal,
if not it should imply $p^2|q-1$ from corollary 2.5 in [Que].
\bn
\item
From   previous lemma \ref{l1304111},  we have
\be\label{e1111101} 
\Big (\frac{x+\zeta^k y}{\mk q_K}\Big)_{\!\! K}
=\Big(\frac{\zeta^{k/2}}{\mk q_K}\Big)_{\!\!
K}\Big(\frac{\epsilon_{k+1}}{\mk q_K}\Big)_{\!\! K} 
\mbox{\ for\ } k=1,\dots,p-2,
\ee
and so, with $p-k$ in place of $k$,
\be\label{e1111102} 
\Big (\frac{x+\zeta^{p-k} y}{\mk q_K}\Big)_{\!\! K}
=\Big(\frac{\zeta^{(p-k)/2}}{\mk q_K}\Big)_{\!\!
K}\Big(\frac{\epsilon_{p-k+1}}{\mk q_K}\Big)_{\!\! K} 
\mbox{\ for\ } p-k=1,\dots,p-2.
\ee
\item  
%We get $\frac{x+\zeta^k y}{x+\zeta^{p-k} y}\in A_k^2\times K^{\timesp}$, so
With the same proof as in  thm \ref{t1107231}, we get 
\be\label{e1111105} 
\Big(\frac{x+\zeta^k y}{\mk q_K}\Big)_{\!\!}
=\Big(\frac{x+\zeta^{p-k}y}{\mk q_K }\Big)_{\!\! K}\mbox{\ for \ }k=2,\dots,p-2,
\ee
which leads  from (\ref{e1111101}) and (\ref{e1111102}) to
\be\label{e1107081}
\Big(\frac{\epsilon_{p-k+1}}{\mk q_K}\Big)_{\!\! K}
=\Big(\frac{\zeta^{k}}{\mk q_K}\Big)_{\!\!
K}\Big(\frac{\epsilon_{k+1}}{\mk q_K}\Big)_{\!\! K} 
\mbox{\ for\ } k=2,\dots,p-2.
\ee
\item
In the other hand,  from (\ref{e1111103}) we  have    
\be\label{e1107082}
 \epsilon_{p-k-1} = \epsilon_{k+1}:
\ee
From (\ref{e1107081}) and (\ref{e1107082}) we derive that
\begin{equation}\label{e1108211}
\Big(\frac{\epsilon_{p-k-1}}{\mk q_K}\Big)_{\!\! K}
=\Big(\frac{\zeta^{-k}}{\mk q_K}\Big)_{\!\!
K}\Big(\frac{\epsilon_{p-k+1}}{\mk q_K}\Big)_{\!\! K}
\mbox{\ for\ k=2,\dots,p-2}.
\ee
\item We get  for the even values $k=2k'$ 
$$\Big(\frac{\epsilon_{p-2k'-1}}{\mk q_K}\Big)_{\!\! K}
=\Big(\frac{\zeta^{-2k'}}{\mk q_K}\Big)_{\!\!
K}\Big(\frac{\epsilon_{p-2k'+1}}{\mk q_K}\Big)_{\!\! K}
\mbox{\ for\ }1\leq k'\leq \frac{p-3}{2}.$$
 Observing that $\epsilon_{p-1}=1$, so 
$\Big(\frac{\epsilon_{p-1}}{\mk q_K}\Big)_{\!\! K}=1$
we get  inductively 
$$\Big(\frac{\epsilon_{p-2k'-1}}{\mk q_K}\Big)_{\!\! K}
=\Big(\frac{\zeta^{-\sum_{j=1}^{k'} 2j}}{\mk q_K}\Big)_{\!\! K}
\Big(\frac{\epsilon_{p-1}}{\mk q_K}\Big)_{\!\! K}
\mbox{\ for\ } k'=1,2,\dots, \frac{p-3}{2},$$
so
$$\Big(\frac{\epsilon_{p-2k'-1}}{\mk q_K}\Big)_{\!\! K}
=\Big(\frac{\zeta^{-k'(k'+1)}}{\mk q_K}\Big)_{\!\! K}
\mbox{\ for\ }0\leq k'\leq \frac{p-3}{2}.$$
\item 
We get for the odd values $k=2k'+1$
$$\Big(\frac{\epsilon_{p-(2k'+1)-1)}}{\mk q_K}\Big)_{\!\! K}
=\Big(\frac{\zeta^{-(2k'+1)}}{\mk q_K}\Big)_{\!\! K}
\Big(\frac{\epsilon_{p-(2k'+1)+1}}{\mk q_K}\Big)_{\!\! K}
\mbox{\ for\ }k' = \frac{p-3}{2},\frac{p-5}{2}\dots,1,$$
so 
$$\Big(\frac{\epsilon_{p-2k'}}{\mk q_K}\Big)_{\!\! K}
=\Big(\frac{\zeta^{2k'+1}}{\mk q_K}\Big)_{\!\! K}
\Big(\frac{\epsilon_{p-2k'-2}}{\mk q_K}\Big)_{\!\! K}
\mbox{\ for\ } k' = \frac{p-3}{2},\frac{p-5}{2}\dots,1.$$
Observing that $\epsilon_1=1$, so 
$\Big(\frac{\epsilon_1}{\mk q_K}\Big)_{\!\! K}=1$
we get for $k'=\frac{p-3}{2}$, so $2k'+1=p-2$,
$$\Big(\frac{\epsilon_{3}}{\mk q_K}\Big)_{\!\! K}
=\Big(\frac{\zeta^{-2}}{\mk q_K}\Big)_{\!\! K}
\Big(\frac{\epsilon_{1}}{\mk q_K}\Big)_{\!\! K},$$
and for $k'=\frac{p-5}{2}$
$$\Big(\frac{\epsilon_{5}}{\mk q_K}\Big)_{\!\! K}
=\Big(\frac{\zeta^{-4}}{\mk q_K}\Big)_{\!\! K}
\Big(\frac{\epsilon_{3}}{\mk q_K}\Big)_{\!\! K},$$
and so on.
\item
Let us define  $k'':=\frac{p-1}{2}-k'$, we get 
$$2k'+1=p-2k'', \mbox{\ for\ } k'=\frac{p-3}{2},\dots,1\mbox{\
corresponding to \ } k''=1, \dots,\frac{p-3}{2}.$$
It follows that  
$$\Big(\frac{\epsilon_{p-2k'}}{\mk q_K}\Big)_{\!\! K}
=\Big(\frac{\zeta^{\sum_{j=1 }^{k''} -2j}}{\mk q_K}
\Big)_{\!\! K}\Big(\frac{\epsilon_1}{\mk q_K}\Big)_{\!\! K}
\mbox{\ for\ } k'= \frac{p-3}{2},\frac{p-5}{2},\dots, 1,$$
so
$$\Big(\frac{\epsilon_{p-2k'}}{\mk q_K}\Big)_{\!\! K}
=\Big(\frac{\zeta^{-k''(k''+1)}}{\mk q_K}
\Big)_{\!\! K}\mbox {\ for\ }1\leq k'\leq \frac{p-3}{2},$$
so
$$\Big(\frac{\epsilon_{p-2k'}}{\mk q_K}\Big)_{\!\! K}
=\Big(\frac{\zeta^{-(\frac{p-1}{2}-k')(\frac{p-1}{2}-k'+1)}}{\mk q_K}
\Big)_{\!\! K}\mbox {\ for\ }1\leq k'\leq \frac{p-3}{2},$$
and finally 
$$\Big(\frac{\epsilon_{p-2k'}}{\mk q_K}\Big)_{\!\! K}
=\Big(\frac{\zeta^{\frac{1}{4}-k'^2}}{\mk q_K}
\Big)_{\!\! K}\mbox {\ for\ }1\leq k'\leq \frac{p-3}{2}.$$
\en
\item
Let us suppose that $q\equiv 1\mod p^2$: then $\Big(\frac{\zeta}{\mk
q_K}\Big)_{\!\! K}=1$ 
and from relation (\ref{e1108211}) we get
$$\Big(\frac{\epsilon_{p-k-1}}{\mk q_K}\Big)_{\!\! K}
=\Big(\frac{\epsilon_{p-k+1}}{\mk q_K}\Big)_{\!\! K}\mbox{\ for\ }
k=2,\dots, p-2.$$
In the other hand we have  $\Big(\frac{\epsilon_{p-1}}{\mk q_K}\Big)_{\!\!
K}
=\Big(\frac{\epsilon_{1}}{\mk q_K}\Big)_{\!\! K}=1$ and so 
$$\Big(\frac{\epsilon_{j}}{\mk q_K}\Big)_{\!\! K}=1\mbox{\  for\
}j=1,\dots,p-1.$$
A straightforward computation shows  that 
$\Big(\frac{\epsilon_1\dots\epsilon_{p-1}}{\mk q_K}\Big)_{\!\! K}
=\Big(\frac{1+\zeta}{\mk q_K}\Big)_{\!\! K}$ and we  
 derive that $$\Big(\frac{1+\zeta}{\mk q_K}\Big)_{\!\! K}=1,$$
and finally that 
$$\Big(\frac{1+\zeta^j}{\mk q_K}\Big)_{\!\! K}=1\mbox{\ for\
}j=1,\dots,p-1.$$
which achieves the proof for $p^2|q-1$.
\en 
\end{proof}
\end{thm}

\paragraph{ Acknowledgments:}
I would like to thank Georges   Gras  for pointing out many errors   
in the preliminary versions and for  suggesting  many  improvements to me
for  the content and form of the article.

Roland Qu\^eme

13 avenue du ch\^ateau d'eau

31490 Brax

France

mailto: roland.queme@gmail.com

\begin{thebibliography}{9}
%\bibitem[BDF] {bel} K. Belabas, F. Diaz Y Diaz, and E. Friedman, {\it 
%Small generators of the ideal class group},
%Mathematics of Computation 77, 262 (2008) 1185--1197.
%\bibitem[Coh]{Coh} H. Cohen,
%\textit{Number Theory Volume 1: Tools and Diophantine equations},
%Springer, 2007.
\bibitem[Fur]{fur} P. Furtw$\ddot{a}$ngler, \textit{Letzter Fermatschen
Satz und Eisensteins'ches Reciprozit$\ddot{a}$tsgesetz}, 
Sitz\-ungsber, Akad. d. Wiss. Wien., Abt. IIa, 121, 1912, 589--592.
\bibitem[Gr1]{gr3} G. Gras, 
\textit{ Class Field Theory, From Theory to Practice}, Springer, 2003.
\bibitem[Gr2]{gr2} G. Gras, 
\textit{ Analysis of the classical cyclotomic approach of Fermat's Last
Theorem},
Publications Math\'ematiques de Besan\c{c}on, 2010.
\bibitem[GQ], G. Gras and R. Qu\^eme, \textit{Vandiver papers on cyclotomy revisited
and Fermat's Last Theorem}, Publications Math\'ematiques de Besan\c con (2012/2), 47-111.
\bibitem [Que], R. Qu\^eme, \textit{On second case of Strong Fermat's Last Theorem conjecture}, preprint submitted arXiv.
\bibitem[Rib1]{rib} P. Ribenboim, 
\textit{13 Lectures on Fermat's Last Theorem}, Springer-Verlag, 1979.
\textit{Introduction to cyclotomic fields, second edition}, Springer, 1997.
\end{thebibliography}
\end{document}